\documentclass[a4paper, 11pt]{article}
\usepackage{amsmath, amssymb, eucal, amscd, amstext}

\newtheorem{thm}{Theorem}[section]

\newtheorem{prop}[thm]{Proposition}
\newtheorem{cor}[thm]{Corollary}
\newtheorem{defn}[thm]{Definition}
\newtheorem{rem}[thm]{Remark}

\newenvironment{prf}{{\noindent \textbf{Proof:}\ }}{\hfill $\Box$\\ \smallskip}

\numberwithin{equation}{section}

\newcommand{\smnoind}{\smallskip\noindent}
\newcommand{\CL}{\mathcal{L}}

\newcommand{\supp}{\ \!{\rm supp\ \!}}

\begin{document}

\title{Property $(T)$ for non-unital $C^*$-algebras\footnote{This work is
jointly supported by Hong Kong RGC Research Grant (2160255),
Taiwan NSC grant (NSC96-2115-M-110-004-MY3) and National Natural Science Foundation of China
(10771106).}}
\author{Chi-Wai Leung, Chi-Keung Ng, Ngai-Ching Wong }

\date{October 9, 2007}

\maketitle

\begin{abstract}
Inspired by the recent work of Bekka, we study two reasonable analogues of
property $(T)$ for not necessarily unital $C^*$-algebras.
The stronger one of the two is called ``property $(T)$'' and the weaker one is called ``property $(T_{e})$''.
It is shown that all non-unital C*-algebras do not have property $(T)$
(neither do their unitalizations).
Moreover, all non-unital $\sigma$-unital C*-algebras do not have property $(T_e)$.
\end{abstract}

\bigskip

\section{Introduction}

\bigskip

Property $(T)$ for locally compact groups was first introduced
and studied by Kazhdan in 1960's (see \cite{Kaz}).
This notion was proved to be very useful in the study of topological groups.
In 1980's, Connes defined the related concept of property $(T)$
for type $\rm II_1$-factors (see \cite{Con}) which was also proved to be important.
Recently, Bekka considered in \cite{B} property $(T)$ for unital
$C^*$-algebras and this was later on studied by Brown in \cite{Br}.

\medskip

In this short article, we will consider property $(T)$ for
not necessarily unital $C^*$-algebras.
Roughly speaking, a unital $C^*$-algebra
$A$ is said to have \emph{property $(T)$} if every Hilbert
$A$-bimodule having an almost central unit vector for $A$ contains a
central unit vector for $A$ (see Section 2).
Notice that in the unital case, in order to check a
C*-algebra having property $(T)$, it suffices to consider only the class
of essential Hilbert bimodules (see Proposition \ref{u>T=wT}(b)).
However, there is no guarantee that it is the case for non-unital
$C^*$-algebras.
For this reason, we introduce two notions of property $(T)$ for
general C*-algebras.
We called the stronger one of the two
``property $(T)$'' and the weaker one ``property $(T_e)$'' (according to whether
we consider the class of all Hilbert bimodules, or we restrict our
attention to essential Hilbert bimodules).

\medskip

The main results of this paper is the (somehow discouraging)
fact that any non-unital $\sigma$-unital $C^*$-algebra does not
even have property $(T_e)$.
Other classes of infinite dimensional
$C^*$-algebras definitely not having property $(T_e)$ include abelian $C^*$-algebras, group
C*-algebras of compact groups and AF-algebras.
However the authors do not know whether all non-unital $C^*$-algebras do not have property $(T_{e})$ (if it is the case, then property $(T)$ and property $(T_e)$ are actually equivalent).

\medskip

On the other hand, the authors are grateful to the referee for showing us how to use a generalisation of Brown's
result in \cite[Theorem 3.4]{Br} to give an elegant proof that non-unital $C^*$-algebras do not have property $(T)$ (see Theorem \ref{nonunital}).

\bigskip

\section{Main results}

\bigskip

\emph{Throughout this article, $A$ will
denote a $C^*$-algebra (not necessarily unital) and
$\widehat A$ is the set of all unitarily equivalence classes
of irreducible representations of $A$.}

\medskip

A \emph{Hilbert bimodule over $A$} is a Hilbert space $H$
together with two commuting $*$-homomorphisms $\rho_{\ell}: A
\rightarrow \mathcal{L}(H)$ and $\rho_{r}: A^{op}
\rightarrow \mathcal{L}(H)$ (where $A^{op}$ is the opposite
algebra of $A$).
For $x\in A$ and $\xi \in H$, we shall write $x\cdot\xi = \rho_{\ell}(x)(\xi)$
and $\xi \cdot x = \rho_{r}(x^{op})(\xi)$.
A Hilbert bimodule $H$ over $A$ is said to be \emph{essential} if the linear
span of $\{x\cdot \xi \cdot y \ |\ x, y \in A\ \mbox{and}\ \xi \in H\}$ is dense
in $H$.
On the other hand, a net $(\xi_{i})$ of unit vectors in $H$
is called \emph{an almost central unit vector} if $\| a\cdot \xi_{i} - \xi_{i}\cdot a\|
\rightarrow 0$ for all $a \in A$.
Moreover, an element $\xi\in H$ is said to be \emph{central}
if $a\cdot \xi = \xi \cdot a$ for all $a\in A$.

\medskip

\begin{defn}
A $C^*$-algebra $A$ is said to have \emph{property $(T)$}
(respectively, \emph{property $(T_e)$}) if every Hilbert bimodule
(respectively, essential Hilbert bimodule) over $A$ having
an almost central unit vector will contain a non-zero central vector.
\end{defn}

\medskip

It is clear that if $A$ has property $(T)$, then $A$ has property $(T_e)$.
In contrary to the unital case (in \cite[Remark 17]{B}), a $C^*$-algebra without tracial
state may not have property $(T_e)$ as can be seen in Theorem \ref{nu+sp>nwT} below.

\medskip

Let us first give the following simple result.

\medskip

\begin{prop}
\label{u>T=wT} Let $A$ be a $C^*$-algebra.

\smnoind (a) $A$ has property $(T)$ if and only if its unitalization $\widetilde{A}$ has
property $(T)$.

\smnoind (b) If $A$ is unital and has property $(T_e)$, then
 $A$ has property $(T)$.
\end{prop}
\begin{prf}
(a) This part is clear.

\smnoind (b) Let $H$ be a Hilbert bimodule over $A$
having an almost central unit vector $(\xi_{i})$.
If $P= \rho_{\ell}(1_{A})$ and $Q = \rho_{r}(1_{A^{op}})$, then:
$$H\ =\ PHQ \oplus PH(1-Q) \oplus (1-P)HQ \oplus (1-P)H(1-Q).$$
It is obvious that $H$ must contain a non-zero central vector if $(1-P)H(1-Q) \neq (0)$, and
so, one can assume that $(1-P)H(1-Q) = (0)$. For each $\xi_{i}$, we write $\xi_{i} =
\alpha_{i} + \beta_{i} +\gamma_{i}$ where $\alpha_{i}\in PHQ$, $\beta_{i}\in PH(1-Q)$ and
$\gamma_{i}\in (1-P)HQ$. Then $$\|x \cdot \alpha_{i} - \alpha_{i} \cdot x\|^{2} + \|x\cdot
\beta_{i}\|^{2} + \|\gamma_{i}\cdot x\|^{2}\ =\ \|x\cdot \xi_{i} -\xi_{i} \cdot x\|^{2} \
\rightarrow\ 0 \qquad (x \in A)$$ will imply that $\beta_{i} \rightarrow 0$ and $\gamma_{i}
\rightarrow 0$ by taking $x=1_{A}$. Since $\|\xi_{i}\| = 1$, we have $\alpha_{i}
\not\rightarrow 0$ and $PHQ$ has an almost central unit vector. As $PHQ$ is an essential
Hilbert bimodule over $A$ and $A$ has property $(T_e)$, we know that $H$ contains a non-zero
invariant vector.
\end{prf}

\medskip

\begin{rem}
{\rm
(a) It is very tempting to use the argument of Proposition \ref{u>T=wT}
(b) to show that property $(T_e)$ is equivalent to property $(T)$.
However, the problem is that if $(a_i)$ is an approximate unit
for $A$, and $P$ and $Q$ are the strong operator limits
of $\rho_l(a_i)$ and $\rho_r(a_i^{op})$ respectively,
it is possible that $PHQ = (0)$,
e.g. $\rho_l: \mathcal{K}(\ell^2) \rightarrow \CL(\ell^2)$ is
the canonical embedding and $\rho_r = 0$ (note that $\{e_n\}$ is an almost central vector).

\smnoind
(b) We do not know if property $(T_e)$ is
also preserved under unitalization. According to Proposition \ref{u>T=wT}(b), this statement
is equivalent to saying that property $(T)$ being the same as property $(T_e)$.}
\end{rem}

\medskip

Part (a) of the following corollary follows from Proposition \ref{u>T=wT} and part (b) is easy to verify.

\medskip

\begin{cor}
\label{sum}
(a) If $B$ and $C$ are $C^*$-algebras having property $(T)$,
then so is their direct sum $B \oplus
C$.

\smnoind
(b) If $A$ has property $(T)$ (respectively, property $(T_e)$),
then so is any quotient of $A$.
\end{cor}

\medskip

Note that if $A$ has property $(T)$, then so is its multiplier algebra $M(A)$. However, the
converse is not true because the $C^*$-algebra of compact operators $\mathcal{K}(\ell^2)$
does not even have property $(T_e)$ (see Proposition \ref{C0-nwT}(a) below) while the $C^*$-algebra of
all bounded linear operators $\mathcal{L}(\ell^2)$ has.

\medskip

In the following, we give several cases when $A$ will definitely not have property $(T_e)$.

\medskip

\begin{prop}
\label{C0-nwT}
Suppose that $A$ is a $C^*$-algebra and $M(A)$ is its
multiplier algebra.

\smnoind
(a) If $\Omega$ is a locally compact Hausdorff
space and $\varphi$ is a non-degenerate $*$-homomorphism
from $C_0(\Omega)$ to $M(A)$ such that $1\notin \varphi(C_0(\Omega))$,
then $A$ does not have property $(T_e)$.

\smnoind
(b) If there exists an infinite directed set $I$ and a net of
increasing projections $\{p_i\}_{i\in I}$ in $M(A)$
such that $p_i \rightarrow 1$ strictly and $p_i \lneq p_j
\lneq 1$ for any $i,j\in I$ with $i\lneq j$, then $A$ does
not have property $(T_e)$.
\end{prop}
\begin{prf}
(a) Let $I := \{K\subseteq \Omega: K$ is compact$\}$.
For any $K\in I$, fix $f_K\in C_c(\Omega)$
with
$$
\chi_K\leq f_K\leq 1_\Omega,
$$
where $\chi_K$ is the characteristic function of $K$. Then by the
non-degeneracy of $\varphi$, one has $\|a - a\varphi(f_K)\| + \|a -
\varphi(f_K)a\| \rightarrow 0$ (along $K\in I$) for any $a\in A$. If
$\pi: A \rightarrow \mathcal{L}(H)$ is any non-degenerate
$*$-representation, it induces a unital $*$-representation
$\varphi_\pi$ of $B(\Omega)$  on $H$ (where $B(\Omega)$ is the
$C^*$-algebra of all bounded Borel measurable functions on
$\Omega$). Suppose that there exists $K\in I$ such that
$\varphi_\pi(\chi_K) = 1$ (and hence $\pi (\varphi(f_K)) = 1$) for
every $\pi\in \widehat A$.
Then the injectivity of $\bigoplus_{\pi\in \widehat A}
\pi$ will imply that $\varphi(f_K) =1$ which contradicts the
hypothesis. This shows that for any $K\in I$, there exists $\pi\in
\widehat A$ such that $\varphi_\pi(\chi_K) \neq 1$.

\smallskip

\noindent \emph{Case 1. There is $\pi\in \widehat A$ such that $\varphi_\pi(\chi_K) \neq 1$ for
any $K\in I$.}

\smallskip

Let $H_\pi$ be the underlying Hilbert space for $\pi$. Assume that $H_\pi$ is finite dimensional.
Then $\pi(\varphi(C_0(\Omega)))$ is unital. It is not hard to check that $\|f -
f\chi_K\|_{B(\Omega)} \rightarrow 0$ for any $f\in C_0(\Omega)$. Thus, $\varphi_\pi(\chi_K)$
converges to $1$ in norm and there exists $K_0\in I$ with $\varphi_\pi(\chi_{K_0}) = 1$
which contradicts the assumption of Case 1. Thus $H_\pi$ is infinite dimensional. It is easy to
check that $HS(H_\pi) \cong H_\pi\otimes \bar H_\pi$ is an essential Hilbert bimodule over $A$ with
multiplications:
$$a\cdot (\xi\otimes \bar \eta)\cdot b\ =\ \pi(a)(\xi) \otimes \overline{\pi(b^*)(\eta)}$$
($\xi, \eta \in H_\pi$). If $\Theta\in HS(H_\pi)$ is a central vector, then $\pi(a) \Theta = \Theta
\pi(a)$ for all $a\in A$.  This implies that $\Theta\in \mathbb{C}  1$ (as $\pi$ is
irreducible) and so $\Theta =0$ (because $H_\pi$ is infinite dimensional). Thus, there is no
non-zero central vector in $HS(H_\pi)$. For any $K\in I$, there exists $\xi_K\in H_\pi$ such that
$\varphi_\pi(\chi_K)\xi_K = 0$ and $\|\xi_K\| =1$.
Define $\zeta_K := \xi_K \otimes
\overline{\xi_K}$.
Then for any $f\in C_0(\Omega)$ with $\supp f\subseteq K$, we have
$\pi(\varphi(f))\xi_K = 0$ and so
\begin{eqnarray*}
\lefteqn{\|a\cdot \zeta_K - \zeta_K \cdot a\|}\\
&\qquad = & \|\pi(a)(\xi_K)\otimes \overline{\xi_K} - \pi(a\varphi(f))(\xi_K)\otimes
\overline{\xi_K} - \xi_K\otimes \overline{\pi(a^*)(\xi_K)} +
\xi_K\otimes \overline{\pi(a^*\varphi(f^*))(\xi_K)}\| \\
& \qquad  \leq & \| a - a\varphi(f)\| + \| a - \varphi(f)a\|.
\end{eqnarray*}
Now for any $\epsilon > 0$, there exists $K_0\in I$ with $\|a - a\varphi(f_{K_0})\| + \|a -
\varphi(f_{K_0})a\| < \epsilon$. For any $K\in I$ with $\supp f_{K_0}\subseteq K$,
$$\|a\cdot \zeta_K - \zeta_K \cdot a\| \leq \|a - a\varphi(f_{K_0})\|
+ \|a - \varphi(f_{K_0})a\| < \epsilon.$$
Consequently, $\{\zeta_K\}_{K\in I}$ is an almost central vector
for $A$ and $A$ does not have property $(T_e)$.

\smallskip

\noindent \emph{Case 2. For any $\pi \in \widehat A$ there exists $K_\pi\in I$ such that
$\varphi_\pi(\chi_{K_\pi}) = 1$.}

\smallskip

Consider $H_0 := \bigoplus_{\underset{\pi,\sigma\in \widehat A}{\pi \neq \sigma}} H_\pi \otimes
\overline{H_\sigma}$ as an essential Hilbert bimodule over $A$ with multiplications: $a\cdot
(\xi\otimes \bar \eta)\cdot b =  \pi(a)(\xi) \otimes \overline{\sigma(b^*)(\eta)}$ for any
$\pi,\sigma\in \widehat A$, $\xi\in H_\pi$ and $\eta\in H_\sigma$. Let $\Theta\in H_0$ be a
central vector. Then $\Theta = (\Theta_{\pi,\sigma})$ where $\Theta_{\pi,\sigma} \in
H_\pi\otimes \overline{H_\sigma}$.
By considering $\Theta_{\pi,\sigma}$ as an element in
$HS(H_\sigma; H_\pi)$, the relation $\pi(a)\Theta_{\pi,\sigma} = \Theta_{\pi,\sigma}
\sigma(a)$ ($a\in A$) and the Schur's lemma tells us that $\Theta_{\pi, \sigma} = 0$ (as
$\pi \neq \sigma$).
This shows that $H_0$ has no non-zero central vector. We claim that for
any $K\in I$, there exist at least two elements $\pi, \sigma\in \widehat A$ such that
$$\varphi_\pi(\chi_K)\neq 1 \qquad {\rm and} \qquad \varphi_\sigma(\chi_K)\neq 1.$$
Indeed, as noted above, there exists $\pi\in \widehat A$ such that
$\varphi_\pi(\chi_K) \neq 1$. Suppose on the contrary that $\varphi_\sigma(\chi_K) = 1$ for
any $\sigma\in \widehat A \setminus \{\pi\}$. Then $\varphi_\sigma(\chi_L) = 1$ for any $L\in I$
with $K \subseteq L$ and any $\sigma\in \widehat A \setminus \{\pi\}$. If $K_\pi$ is as in the
assumption of Case 2 and if $L\in I$ with $K\subseteq L$ and $K_\pi\subseteq L$, then
$\bigoplus_{\sigma\in \widehat A} \sigma(\varphi(f_L)) = 1$ which
contradicts the hypothesis that $1\notin \varphi(C_0(\Omega))$. Now for any $K\in I$, we
take two different elements $\pi,\sigma \in \widehat A$ with $\varphi_\pi(\chi_K)\neq 1$ and
$\varphi_\sigma (\chi_K) \neq 1$, and we choose $\xi_K\in H_\pi$ and $\eta_K\in H_\sigma$
such that $\varphi_\pi(\chi_K)(\xi_K) = 0$,
$\varphi_\sigma(\chi_K)(\eta_K) = 0$ and $\|\xi_K\| = 1 = \|\eta_K\|$.
Define $\zeta_K := \xi_K\otimes \overline{\eta_K}$.
For any
$f\in C_0(\Omega)$ with $\supp f\subseteq K$, we have
\begin{eqnarray*}
\lefteqn{\|a\cdot \zeta_K - \zeta_K \cdot a\|} \\
& \qquad = &\|\pi(a)(\xi_K)\otimes \overline{\eta_K} - \pi(a\varphi(f))(\xi_K)\otimes
\overline{\eta_K} - \xi_K\otimes \overline{\sigma(a^*)(\eta_K)} +
\xi_K\otimes \overline{\sigma(a^*\varphi(f^*))(\eta_K)}\| \\
& \qquad \leq & \| a - a\varphi(f)\| + \| a - \varphi(f)a\|.
\end{eqnarray*}
Now a similar argument as that of Case 1 will show that $A$ does not have property $(T_e)$.

\medskip
\noindent
(b) The proof is similar to that of part (a) (but we need to replace $\{\chi_K\}$ with $\{p_i\}$).
\end{prf}

\medskip

Although property $(T)$ is preserved under a finite direct sum (see Corollary \ref{sum}), it
does not hold for an infinite $c_{0}$-direct sum. More precisely, we have the following
direct application of Proposition \ref{C0-nwT}(b).

\medskip

\begin{cor}
\label{c0} If $(A_{\lambda})_{\lambda\in \Lambda}$ is any infinite family of nonzero
$C^*$-algebras, then the $c_0$-direct sum $\bigoplus_{\lambda\in \Lambda} A_{\lambda} :=
\left\{ (x_{\lambda})_{\lambda\in \Lambda}\in \Pi_{\lambda\in \Lambda} A_{\lambda}:
(\|x_{\lambda}\|)_{\lambda\in \Lambda} \in c_0(\Lambda)\right\}$ does
not have property $(T_{e})$.
\end{cor}

\medskip

Suppose that a non-untial $C^*$-algebra $A$ contains a strictly positive
element $h$ (see \cite[3.10.6]{Ped}).
The smallest $C^*$-subalgebra $B\subseteq A$ generated by $h$
is isomorphic to $C_0(\Omega)$ for some non-compact locally compact space.
Since $\{h^{1/n}\}$ is an approximate identity for $A$,
Proposition \ref{C0-nwT}(a) gives the following result (which implies that property $(T_{e})$ is equivalent to property $(T)$ for any $\sigma$-unital $C^*$-algebra).

\medskip

\begin{thm}
\label{nu+sp>nwT} Every non-unital $\sigma$-unital $C^*$-algebra (in particular, any separable
non-unital $C^*$-algebra) does not have property $(T_{e})$. 
\end{thm}

\medskip

Proposition \ref{C0-nwT} also gives the following corollary.
Part (a) of it follows from Proposition
\ref{C0-nwT}(a) and \cite[Proposition 15]{B} while part (b) follows from \cite[Theorem
28.40]{HR} and Corollary \ref{c0}. To show part (c), one needs (on top of Theorem
\ref{nu+sp>nwT}) \cite[Proposition 5.1]{Br}
 as well as the fact that any unital AF-algebra has a tracial state (see \cite{Haag}).

\medskip

\begin{cor}
Let $A$ be a $C^*$-algebra.
If $A$ is in one of the following three classes of
$C^*$-algebras, then $A$ having property $(T_e)$ will imply that
$A$ is finite dimensional: (a) $A$ is commutative; (b) $A  = C^*(G)$ for a compact group
$G$; (c) $A$ is an AF-algebra.
\end{cor}

\medskip


It is believed that one can remove the $\sigma$-unital assumption in
Theorem \ref{nu+sp>nwT} (note that the two cases considered in
Proposition \ref{C0-nwT} do not have such assumption) although we still do not have a proof.
However, if only property $(T)$, instead of property $(T_e)$, is concerned, the referee has kindly informed us that this is true (see the following theorem).
As an application, we see that if $A$ is a non-unital $C^*$-algebra,
then $\widetilde A$ will never have property $(T)$.

\medskip

\begin{thm}\label{nonunital}
All non-untial $C^*$-algebras do not have property $(T)$.
\end{thm}
\begin{prf}[Provided by the referee.]
In \cite[Theorem 3.4]{Br}, N.P. Brown showed that if $B$ is a
separable unital $C^*$-algebra with property $(T)$ and $\pi: B\rightarrow M_{n}({\mathbb
C})$ is any irreducible representation, then the central cover $c(\pi)$
(i.e., a central projection in $B^{**}$ defined by
$B^{**}(1-c(\pi)) = \ker \pi^{**}$) of $\pi$ must belong to $B$. Indeed, the separability
assumption  can be removed by replacing Voiculescu's Theorem with
Glimm's Lemma \cite[Lemma II.5.1]{D} in the proof of \cite[Theorem 3.4]{Br}.

\smallskip
Let $A$ be a non-unital $C^*$-algebra.  Suppose on the contrary that
$A$ had property $(T)$. Then its unitization $\widetilde{A}$ also has
property $(T)$.
Let $\pi : \widetilde{A} \rightarrow {\mathbb C}$ be the canonical map.
By the extension of \cite[Theorem 3.4]{Br} as stated above, the central cover $c(\pi)$ of $\pi$ is contained in
$\widetilde{A}$. This yields the following $C^*$-algebras decompositions:
$$
\widetilde{A} =\ (1-c(\pi))\widetilde{A} \oplus c(\pi)\widetilde{A} = \ker \pi
\oplus {\mathbb C} = A \oplus {\mathbb C},
$$
which implies the contradiction that  $A$ is unital.
\end{prf}

\medskip\noindent
Chi-Wai Leung, Department of Mathematics, The Chinese University of Hong Kong, Hong
Kong.

\smnoind \emph{Email address:} cwleung@math.cuhk.edu.hk

\medskip\noindent
Chi-Keung Ng, Chern Institute of Mathematics and LPMC, Nankai University, Tianjin 300071,
China.

\smnoind \emph{Email address:} ckng@nankai.edu.cn

\medskip\noindent
Ngai-Ching Wong, Department of Applied Mathematics, National Sun Yat-sen
University, and National Center for Theoretical Sciences, National Science Council,
Kaohsiung, 804, Taiwan. (Current address:
Department of Mathematics, The Chinese University of Hong Kong, Shatin,
New Territories, Hong Kong.)

\smnoind \emph{Email address:} wong@math.nsysu.edu.tw

\end{document}